\documentclass[a4, 12pt]{amsart}
\usepackage[mathscr]{eucal}
\usepackage{amssymb}
\usepackage{latexsym}
\usepackage{amsthm}
\theoremstyle{plain}
\newtheorem{theorem}{Theorem}[section]
\newtheorem{corollary}{Corollary}[section]
\newtheorem{remark}{Remark}[section]
\newtheorem{Lemma}{Lemma}[section]

\numberwithin{equation}{section} 
\makeatletter
\title[eigenvalues  of the Paneitz operator]{Estimates for  eigenvalues  of  \\ the Paneitz operator*}
\author{Qing-Ming Cheng} 
\address{\noindent Department of Applied Mathematics, Faculty of Sciences, 
\newline \indent Fukuoka University, Fukuoka  814-0180, Japan \newline \indent cheng@fukuoka-u.ac.jp}

\begin{document}
\maketitle

\begin{abstract}
\noindent For an  $n$-dimensional compact 
submanifold $M^n$ in the Euclidean space $\mathbf R^{N}$, we study estimates for  eigenvalues of
the Paneitz operator on $M^n$. Our  estimates for eigenvalues are sharp. 
\end{abstract}



\footnotetext{* Research partially Supported by JSPS Grant-in-Aid for
Scientific Research (B) No. 24340013.}

\section {Introduction}

\noindent
For   compact  Riemann surfaces $M^2$, Li and Yau \cite{LY}
introduced the notion of conformal volume, which is a global invariant of the conformal structure. 
They determined the conformal volume for a large class of Riemann surfaces, which
admit minimal immersions into spheres. In particular, they proved that for a compact Riemann surface $M^2$, if there exists
a conformal map from $M^2$ into the unit sphere $S^N(1)$, then the first eigenvalue $\lambda_1$ of the Laplacian
satisfies
$$
\lambda_1\text{vol}(M^2)\leq 2V_c(N,M^2)
$$
and the equality holds only if $M^2$ is a minimal surface in  $S^N(1)$, 
where $V_c(N,M^2)$ is the conformal volume of $M^2$.

\noindent
For $4$-dimensional compact 
Riemannian manifolds, Paneitz \cite{P} introduced a fourth order operator $P_g$
defined by, letting $\text{div}$ be the divergence for the metric $g$,  
\begin{equation}
P_gf = \Delta^2f -\text{div}\bigl [(\frac23Rg-2\text{Ric})\nabla f\bigl],
\end{equation}
for smooth functions $f$ on $M^4$, where  $\Delta$ and $\nabla$ denote  the Laplacian  
and the gradient operator  with respect to the metric $g$ on $M^4$, respectively,  and 
$R$ and $\text{Ric}$ are the scalar curvature and Ricci curvature tensor with respect to the metric $g$
on $M^4$. 
Furthermore, Branson \cite{B} has generalized the Paneitz operator to an $n$-dimensional 
Riemannian manifold. For an $n$-dimensional Riemannian manifold $(M^n,g)$, the operator
$P_g$ is defined by 
\begin{equation}
P_gf = \Delta^2f -\text{div}\bigl [(a_nRg+b_n\text{Ric})\nabla f\bigl] +\dfrac{n-4}2Qf,
\end{equation}
where 
$$
Q=c_n|\text{Ric}|^2+d_n-\dfrac{1}{2(n-1)}\Delta R
$$
is called $Q$-curvature with respect to the metric $g$,
$$
a_n=\dfrac{(n-2)^2+4}{2(n-1)(n-2)}, \ \ b_n=-\frac{4}{n-2}, \
$$
$$  c_n=-\dfrac{2}{(n-2)^2}, \  \ 
d_n=\dfrac{n(n-2)^2-16}{8(n-1)^2(n-2)^2}.
$$
This operator $P_g$ is also called Paneitz operator or Branson-Paneitz operator.
It is known that  Paneitz operator is conformally invariant of bi-degree $(\dfrac{n-4}2, \dfrac{n+4}2)$, that is, under conformal 
transformation
of Riemannian metric  $g=e^{2w}g_0$, the Paneitz operator $P_g$ changes into
\begin{equation}
P_gf=e^{-\frac{n+4}2w}P_{g_0}(e^{\frac{n-4}2w}f).
\end{equation}

\noindent
Let $\frak{M}(M^n)$ be the set of Riemannian metrics on $M^n$. For each $ g \in \frak{M}(M^n)$,
the total $Q$-curvature for $g$ is defined by
\begin{equation*}
Q[g]=\int_{M^n}Qdv.
\end{equation*}
When $n=4$, from the Gauss-Bonnet theorem for dimension $4$, we have 
\begin{equation}
Q[g]=-\frac14\int_{M^4}|W|^2dv+8\pi^2\chi(M^4),
\end{equation}
where $W$ is the Weyl conformal curvature tensor and $\chi(M^4)$ is the Euler characteristic of $M^4$.
Hence, we know that the total $Q$-curvature is a conformal invariant for dimension $4$. In \cite{N}, Nishikawa has 
studied the variation of the total $Q$-curvature for a general dimension $n$. He has proved that a Riemannian metric $g$ 
on an $n (n\neq 4)$-dimensional compact manifold $M^n$ is a critical point of the total $Q$-curvature functional with respect 
to  a volume preserving conformal variation of the  metric $g$, if and only if the $Q$-curvature with respect to the metric 
$g$ is constant. 

\noindent
Since  the Paneitz operator  $P_g$ is an elliptic operator and $P_g1=0$ for $n=4$, we know that $\lambda_0=0$ is 
an eigenvalue of $P_g$. Gursky \cite{G1} shown that if the Yamabe invariant of $M^4$ is positive and the total $Q$-curvature
is positive, the first eigenvalue $\lambda_1$ is positive.  For $n\geq 6$, Yang and Xu \cite{XY1} have proved the 
Paneitz operator $P_g$ is positive if the scalar curvature is positive and $Q$-curvature is nonnegative. Furthermore, 
see \cite{BCY, CY, CGY, HR}.

\noindent
For $n\geq 3$, we consider the following closed eigenvalue problem on an $n$-dimensional compact
manifold $M^n$:
\begin{equation}
P_gu = \lambda u.
\end{equation}
Since $P_g$ is an elliptic operator, the spectrum of $P_g$ on $M^n$ is discrete. We assume 
$$
0<\lambda_1<\lambda_2\leq \cdots, \lambda_k\leq \cdots\to +\infty
$$
for $n\neq 4$ and for $n=4$, 
$$
0=\lambda_0<\lambda_1\leq \lambda_2\leq \cdots, \lambda_k\leq \cdots\to +\infty.
$$
When $n=4$, Yang and Xu \cite{XY2} have introduced an $N$-conformal energy $E_c(N,M^4)$ if $M^4$ can
be conformally immersed into the unit sphere $S^N(1)$ and have obtained  an upper bound for the first eigenvalue
$\lambda_1$:
$$
\lambda_1\text{vol}(M^4) \leq E_c(N,M^4),
$$
where $\text{vol}(M^n)$ denotes the volume of $M^n$. Furthermore, Chen and Li \cite{CL1} have also studied
the upper bound on the first eigenvalue
$\lambda_1$ when $M^4$ is considered as a compact submanifold in a Euclidean space $\mathbf R^N$.  They have 
proved
\begin{equation*}
 \aligned &
 \lambda_{1}\leq \dfrac{\int_{M^4}\bigl(16|H|^2+\frac23R\bigl)dv
 \int_{M^4}|H|^2dv}{\bigl\{\text{vol}(M^4)\bigl\}^2}\\
 \endaligned
\end{equation*}
and the equality holds if and only if   $M^4$ is a minimal submanifold in a sphere $S^{N-1}(r)$ for $N>5$
and  $M^4$ is a round  sphere $S^{4}(r)$ for $N=5$. 
In \cite{CL2}, the second eigenvalue $\lambda_2$ of the Paneitz operator $P_g$ is studied. By making use of 
the conformal transformation introduced by Li and Yau \cite{LY}, Chen and Li  proved, for $n\geq 7$,
$$
\lambda_2\text{vol}(M^n)\leq \dfrac12n(n^2-4)\int_{M^n}\bigl |H|^4dv+\dfrac{n-4}2\int_{M^n}Qdv
$$
 if  $M^n$ is a compact submanifold in the Euclidean space $\mathbf R^N$. Here $|H|$ denotes the mean curvature of $M^n$ in $\mathbf R^N$. As they remarked,
 their method does not work for $3\leq n\leq 6$. 
 
 \noindent
 The purpose of this paper is to study eigenvalues of the Paneitz operator $P_g$ in  $n$-dimensional compact 
 Riemannian manifolds. Our method is very different from one used by Chen and Li \cite{CL2} 
 and Xu and Yang \cite{XY2}. From Nash's theorem, we know that  each compact  
 Riemannian manifold can be isometrically immersed into a Euclidean space $\mathbf {R}^{N}$. 
 Thus, we can assume $M^n$ is  an $n$-dimensional compact  submanifold in $\mathbf {R}^{N}$.
\begin{theorem}
Let $(M^4, g)$ be a  $4$-dimensional compact  submanifold with the 
 metric $g$ in $\mathbf {R}^{N}$. Then,  eigenvalues
of the Paneitz operator $P_g$ satisfy
\begin{equation*}
\aligned &
\sum\limits_{j=1}^4\lambda_{j}^{\frac{1}{2}}\leq 4\dfrac{\sqrt{\int_{M^4}\bigl(16|H|^2+\frac23R\bigl)dv
 \int_{M^4}|H|^2dv}}{\text{\rm vol}(M^4)}\\
 \endaligned
\end{equation*}
and the equality holds if and only if 
 $M^4$ is a  round sphere  $S^{4}(r)$ for  $N=5$ and  $M^4$ is a  compact minimal submanifold with constant
scalar  curvature in  $S^{N-1}(r)$ for  $N>5$. 
\end{theorem}

\begin{corollary}
Let $(M^4, g)$ be a  $4$-dimensional compact  submanifold with the 
 metric $g$ in the unit sphere  $S^{N}(1)$. Then,  eigenvalues
 of the Paneitz operator $P_g$ satisfy
\begin{equation*}
\aligned &
\sum\limits_{j=1}^4\lambda_{j}^{\frac{1}{2}}\leq 4\dfrac{\sqrt{\int_{M^4}\bigl(16|H|^2+16+\frac23R\bigl)dv
 \int_{M^4}(|H|^2+1)dv}}{\text{\rm vol}(M^4)}\\
 \endaligned
\end{equation*}
and the equality holds if and only if 
 $M^4$ is a compact minimal submanifold with constant scalar curvature in   $S^{N}(1)$.
\end{corollary}
\begin{theorem}
Let $(M^n, g)$ $(n>4)$ be an  $n$-dimensional compact  submanifold with the 
 metric $g$ in  $\mathbf R^{N}$. Then, eigenvalues 
of the Paneitz operator $P_g$ satisfy
 \begin{equation*}
 \aligned 
  &\sum\limits_{j=1}^n(\lambda_{j+1}-\lambda_1)^{\frac{1}{2}}\\
  &\leq\sqrt{\int_{M^n}\dfrac{n(n^2-4)|H|^2}2u_1^2dv+2(n+2)\int_{M^n}g( \nabla u_1, \nabla  u_1)dv}\\
&\times\sqrt{\int_{M^n}n^2|H|^2u_1^2dv+4\int_{M^n}g( \nabla u_1, \nabla  u_1)dv}
\endaligned
\end{equation*}
and  the equality holds if and only if  $M^n$ is isometric to a sphere $S^n(r)$,
where $u_1$ is the normalized first eigenfunction of $P_g$.
\end{theorem}

\begin{remark} In our theorem 1.2, we do not need to assume the positivity of the Paneitz operaator $P_g$.
\end{remark}

\noindent
If the Paneitz operator $P_g$ is a positive operator, we have 
\begin{theorem}
Let $(M^n, g)$ $(n\neq 4)$ be an  $n$-dimensional compact  submanifold with the 
 metric $g$ in the unit sphere  $S^{N}(1)$. Then, eigenvalues 
of the Paneitz operator $P_g$ satisfy
\begin{equation*}
 \aligned &
 \sum_{j=1}^n\lambda_{j}^{\frac{1}{2}}<n\dfrac{\sqrt{\int_{M^n}\bigl((n^2|H|^2+n^2)+(na_n+b_n)R
 +\dfrac{n-4}2Q\bigl)dv\int_{M^n}(|H|^2+1)dv}}{\text{\rm vol}(M^n)}.\\
 \endaligned
\end{equation*}

\end{theorem}

\vskip 10mm
\section{Eigenvalues of the Paneitz operator  on $M^4$}
\noindent
Since  $M^n$ is an $n$-dimensional  submanifold in $\mathbf{R}^N$.
Let ($x_1, \cdots, x_n$) be a local coordinate system in a neighborhood $U$ of $p\in M^n$. Let
${\bf y}$ be the position vector of $p$ in $\mathbf{R}^N$,  which is defined by
$$
{\bf y}=(y_1(x_1, \cdots, x_n),\cdots, y_N(x_1, \cdots, x_n)).
$$
Let  $g$ denote the induced metric of $M^n$ from $\mathbf{R}^N$ and $< , >$ is the standard inner product in
$\mathbf{R}^N$. Thus, we have

\begin{Lemma} {\it For any function $u\in
C^{\infty}(M^n)$, we have
\begin{equation}
\begin{aligned}
&g_{ij}=g(\frac{\partial}{\partial x_i}, \frac{\partial}{\partial
x_j})=<\sum\limits_{\alpha=1}^N\frac{\partial y_{\alpha}}{\partial
x_i}\frac{\partial}{\partial y_{\alpha}},
\sum\limits_{\beta=1}^N\frac{\partial y_{\beta}}{\partial
x^i}\frac{\partial}{\partial
y_{\beta}}>=\sum\limits_{\alpha=1}^N\frac{\partial
y_{\alpha}}{\partial x^i}\frac{\partial y_{\alpha}}{\partial
x^j},\\
&\sum\limits_{\alpha=1}^N(g(\nabla y_{\alpha}, \nabla u))^2=|\nabla
u|^2, \\
&\sum\limits_{\alpha=1}^Ng(\nabla y_{\alpha}, \nabla
y_{\alpha})=\sum\limits_{\alpha=1}^N|\nabla
y_{\alpha}|^2=n,\\
&\sum\limits_{\alpha=1}^N(\Delta
y_{\alpha})^2=n^2|H|^2,\\
&\sum\limits_{\alpha=1}^N\Delta y_{\alpha}\nabla
y_{\alpha}=0,\\
\end{aligned}
\end{equation}
where $\nabla$ denotes the gradient
operator on $M^n$ and $|H|$ is the mean curvature of $M^n$. }
\end{Lemma}

\vskip 3pt\noindent {\it Proof of  Theorem 1}. 
 Let
$u_i$ be an eigenfunction corresponding to  eigenvalue
$\lambda_i$ such that $\{u_i\}_{i=0}^{\infty}$ becomes an
orthonormal basis of $L^2(M^n)$, that is,
\begin{equation*}
\begin{cases}
P_gu_i=\lambda_i u_i, \\
\int_{M^4}u_iu_jdv=\delta_{ij}, \ \ i, j= 0, 1, \cdots.
\end{cases}
\end{equation*}
We define an $N\times N$-matrix $A$ as follows:
$$
A:=(a_{\alpha\beta})
$$
where $a_{\alpha\beta}=\int_{M^4}y_{\alpha}u_0u_{\beta}dv$, for $\alpha, \beta =1, 2, \cdots, N$,  and
${\bf y}=(y_{\alpha})$ is the position vector of the immersion in
$\mathbf{R}^N$. Using the orthogonalization of Gram and Schmidt, we
know that there exist an upper triangle matrix 
$T=(T_{\alpha\beta})$
and an orthogonal matrix 
$U=(q_{\alpha\beta})$ such that $T=UA$,
i.e.,
$$
T_{\alpha\beta}=\sum\limits_{\gamma=1}^Nq_{\alpha\gamma}a_{\gamma\beta}=
\int_{M^4}\sum\limits_{\gamma=1}^Nq_{\alpha\gamma}y_{\gamma}u_0u_{\beta}dv=0,\
\ 1\leq\beta<\alpha\leq N.
$$
Defining
$z_{\alpha}=\sum\limits_{\gamma=1}^Nq_{\alpha\gamma}y_{\gamma}$, we
get
$$
\int_{M^4}z_{\alpha}u_0u_{\beta}dv=\int_{M^4}\sum\limits_{\gamma=1}^Nq_{\alpha\gamma}
y_{\gamma}u_0u_{\beta}dv=0,\ \ 1\leq\beta<\alpha\leq
N.
$$
Putting 
$$
\psi_{\alpha}:=(z_{\alpha}-b_{\alpha})u_0,\ \ \  \ \
b_{\alpha}:=\int_{M4}z_{\alpha}u_0^2dv,\ \ \ \ \ \ 1\leq
\alpha\leq N,
$$ 
we infer 
$$\int_{M^4}\psi_{\alpha}u_{\beta}dv=0,\ \ \ \
0\leq\beta<\alpha\leq N.
$$
Thus, 
from the Rayleigh-Ritz inequality, we have
$$
\lambda_{\alpha}{\int_{M^4}\psi_{\alpha}^2}dv\leq{\int_{M^4}\psi_{\alpha}P_g\psi_{\alpha}}dv,\
\ 1\leq\alpha\leq N.
$$
Since $u_0$ is constant and 
\begin{equation}
\aligned
P_g\psi_{\alpha}&= \Delta^2(z_{\alpha}u_0) -\text{div}\bigl [(\frac23Rg-2\text{Ric})\nabla (z_{\alpha}u_0)
\bigl], 
\endaligned
\end{equation}
according to the Stokes formula, we derive
\begin{equation*}
\aligned
\int_{M^4}\psi_{\alpha}P_g\psi_{\alpha}dv&=\int_{M^4}\biggl[\bigl(\Delta z_{\alpha}\bigl)^2u_0^2
+g((\frac23Rg-2\text{Ric})\nabla z_{\alpha}, \nabla z_{\alpha})u_0^2\biggl]dv.\\
       \endaligned
\end{equation*}
From the lemma 2.1, we have
\begin{equation*}
\aligned
&\sum_{\alpha=1}^N\int_{M^4}\psi_{\alpha}P_g\psi_{\alpha}dv\\&
=\sum_{\alpha=1}^N\int_{M^4}\biggl[\bigl(\Delta z_{\alpha}\bigl)^2u_0^2
+g((\frac23Rg-2\text{Ric})\nabla z_{\alpha}, \nabla z_{\alpha})u_0^2\biggl]dv\\
&=\int_{M^4}\bigl(16|H|^2+\frac23R\bigl)u_0^2dv.
       \endaligned
\end{equation*}
Hence, 
\begin{equation}
\aligned
\sum_{\alpha=1}^N\lambda_{\alpha}{\int_{M^4}\psi_{\alpha}^2}dv
\leq \int_{M^4}\bigl(16|H|^2+\frac23R\bigl)u_0^2dv. 
\endaligned
\end{equation}
On the other hand, 
\begin{equation*}
\aligned &\int_{M^4}\psi_{\alpha}(u_0\Delta z_{\alpha})dv\\
          &=\int_{M^4}( z_{\alpha}u_0-u_0b_{\alpha})(u_0\Delta z_{\alpha})dv\\
                     &=-\int_{M^4}|\nabla(z_{\alpha}u_0)|^2dv.
\endaligned
\end{equation*}
Therefore, for any positive $\delta>0$, we obtain from (2.3)
\begin{equation*}
 \aligned &
 \lambda_{\alpha}^{\frac12}\int_{M^4}|\nabla(z_{\alpha}u_0)|^2dv\\
&=-\lambda_{\alpha}^{\frac12}\int_{M^4}\psi_{\alpha}(u_0\Delta z_{\alpha})dv\\
 &\leq \frac12\bigl(\delta\lambda_{\alpha}\int_{M^4}\psi_{\alpha}^2dv
 +\dfrac{1}{\delta} \int_{M^4}(u_0\Delta z_{\alpha})^2dv\bigl)\\
 \endaligned
\end{equation*}

\begin{equation}
 \aligned &
 \sum_{\alpha=1}^N\lambda_{\alpha}^{\frac12}\int_{M^4}|\nabla(z_{\alpha}u_0)|^2dv\\
 &\leq \frac12\bigl(\delta\sum_{\alpha=1}^N\lambda_{\alpha}\int_{M^4}\psi_{\alpha}^2dv
 +\dfrac{1}{\delta}\sum_{\alpha=1}^N \int_{M^4}(u_0\Delta z_{\alpha})^2dv\bigl)\\
 &\leq \frac12\bigl(\delta\int_{M^4}\bigl(16|H|^2+\frac23R\bigl)u_0^2dv
 +\dfrac{1}{\delta}\int_{M^4}16|H|^2u_0^2dv\bigl).\\
 \endaligned
\end{equation}
It is not hard to prove that,  for any point and  for any $\alpha$,
$$
\aligned
|\nabla z_{\alpha}|^2&=g(\nabla z_{\alpha}, \nabla z_{\alpha})\leq 1.
\endaligned
$$
Hence, 
\begin{equation}
\aligned
&\sum_{\alpha=1}^N\lambda_{\alpha}^{\frac{1}{2}}|\nabla z_{\alpha}|^2\\
   &\geq\sum_{i=1}^4\lambda_{i}^{\frac{1}{2}}|\nabla
   z_{i}|^2+\lambda_{5}^{\frac{1}{2}}\sum\limits_{A=5}^N|\nabla z_{A}|^2\\
   &=\sum_{i=1}^4\lambda_{i}^{\frac{1}{2}}|\nabla
   z_{i}|^2+\lambda_{5}^{\frac{1}{2}}(4-\sum\limits_{j=1}^4|\nabla
   z_{j}|^2)\\
   &\geq \sum_{i=1}^4\lambda_{i}^{\frac{1}{2}}|\nabla z_{i}|^2+
   \sum\limits_{j=1}^4\lambda_{j}^{\frac{1}{2}}(1-|\nabla z_{j}|^2)
   \\
     &\geq \sum\limits_{j=1}^4\lambda_{j}^{\frac{1}{2}}.
\endaligned
\end{equation}
We obtain, by (2.4) and (2.5), 
\begin{equation*}
 \aligned &
  \sum\limits_{j=1}^4\lambda_{j}^{\frac{1}{2}} &\leq \frac12\bigl(\delta\int_{M^4}\bigl(16|H|^2+\frac23R\bigl)u_0^2dv
 +\dfrac{1}{\delta}\int_{M^4}16|H|^2u_0^2dv\bigl).\\
 \endaligned
\end{equation*}
Taking 
$$
\dfrac1{\delta} =\sqrt{\dfrac{\int_{M^4}\bigl(16|H|^2+\frac23R\bigl)u_0^2dv}
{ \int_{M^4}16|H|^2u_0^2dv}}
 $$
 we have, because of  $u_0=\sqrt{\frac1{\text{vol}(M^4)}}$,
\begin{equation}
 \aligned &
\sum\limits_{j=1}^4\lambda_{j}^{\frac{1}{2}}\leq 4\dfrac{\sqrt{\int_{M^4}\bigl(16|H|^2+\frac23R\bigl)dv
 \int_{M^4}|H|^2dv}}{\text{vol}(M^4)}.\\
 \endaligned
\end{equation}
If the equality holds,  we have 
$$
\lambda_1=\lambda_2=\cdots=\lambda_{N},
$$
\begin{equation}
\Delta (z_{\alpha}-b_{\alpha})=-\sqrt{\lambda_5} \delta (z_{\alpha}-b_{\alpha}).
\end{equation}
According to Takahashi's theorem, we know that   $M^4$ is a round  sphere $S^{4}(r)$ for $N=5$
and  $M^4$ is a minimal submanifold in a sphere $S^{N-1}(r)$ for $N>5$ with 
$\sum_{\alpha=1}^N(z_{\alpha}-b_{\alpha})^2=r^2$. 
Thus, we have 
$$
\lambda_1=\lambda_2=\cdots=\lambda_N=\dfrac{16}{r^4\delta^2}.
$$
From the definition of the Paneitz operator $P_g$, we have 
\begin{equation}
P_g(z_{\alpha}-b_{\alpha}) = \Delta^2(z_{\alpha}-b_{\alpha}) -\text{div}
\bigl [(\frac23Rg-2\text{Ric})\nabla(z_{\alpha}-b_{\alpha})\bigl],
\end{equation}
that is, from (2.7) and (2.8), we have 
\begin{equation*}
\lambda_5(1-\delta^2)(z_{\alpha}-b_{\alpha})= -\text{div}
\bigl [(\frac23Rg-2\text{Ric})\nabla(z_{\alpha}-b_{\alpha})\bigl].
\end{equation*}
According to $\sum_{\alpha=1}^N(z_{\alpha}-b_{\alpha})^2=r^2$, we obtain
\begin{equation*}
\lambda_5(1-\delta^2)r^2
=\sum_{\alpha=1}^Ng((\frac23Rg-2\text{Ric})\nabla(z_{\alpha}-b_{\alpha}), \nabla(z_{\alpha}-b_{\alpha})).
\end{equation*}
Hence,
\begin{equation*}
\lambda_5(1-\delta^2)r^2
=\frac23R.
\end{equation*}
Thus, the scalar curvature  $R$ is constant. Hence, $M^4$ is a compact minimal submanifold with constant scalar curvature
in a sphere $S^{N-1}(r)$. 
 This finishes the proof of theorem 1.1.

\vskip .5cm
\noindent
{\it Proof of Corollary 1.1.} Since the unit sphere $S^N(1)$ is a hypersurface 
in $\mathbf R^{N+1}$ with the mean curvature $1$, $M^4$ can be seen as 
a compact submanifold in $\mathbf R^{N+1}$ with the mean curvature $\sqrt{|H|^2+1}$.
According to the theorem 1.1, we complete the proof of the corollary 1.1.

\vskip 1cm
\section{Eigenvalues of the Paneitz operator  on $M^n \ (n\neq 4)$}

\vskip 1cm
\noindent
{\it Proof of theorem 1.2.}
Since $n> 4$, eigenvalues of the Paneitz operator $P_g$
satisfy
$$
\lambda_1< \lambda_2\leq \cdots, \lambda_k\leq \cdots\to +\infty.
$$
Let $u_i$ be an eigenfunction corresponding to  eigenvalue
$\lambda_i$ such that $\{u_i\}_{i=1}^{\infty}$ becomes an
orthonormal basis of $L^2(M^n)$, that is,
\begin{equation*}
\begin{cases}
P_gu_i=\lambda_i u_i, \\
\int_{M^n}u_iu_jdv=\delta_{ij}, \ \ i, j= 1, 2,  \cdots.
\end{cases}
\end{equation*}
We shall use the same idea to prove the theorem 1.2. 
But, in this case, we need to use  the first eigenfunction $u_1$,
which is not constant in general.  Thus, we need to compute many formulas.
We define an $N\times N$-matrix $A$ as follows:
$$
A:=(a_{\alpha\beta})
$$
where $a_{\alpha\beta}=\int_{M^n}y_{\alpha}u_1u_{\beta+1}dv$, for $\alpha, \beta =1, 2, \cdots, N$,  and
${\bf y}=(y_{\alpha})$ is the position vector of the immersion in
$\mathbf{R}^N$. 
Thus, there is  an orthogonal matrix  $U=(q_{\alpha\beta})$ such that 
$$
\int_{M^n}z_{\alpha}u_1u_{\beta+1}dv=0,\ \ 1\leq\beta <\alpha\leq N,
$$
where 
$z_{\alpha}=\sum\limits_{\gamma=1}^Nq_{\alpha\gamma}y_{\gamma}$.
Putting 
$$
\varphi_{\alpha}:=(z_{\alpha}-a_{\alpha})u_1,\ \ \  \ \
a_{\alpha}:=\int_{M^n}z_{\alpha}u_1^2dv,\ \ \ \ \ \ 1\leq
\alpha\leq N,
$$ 
we infer 
$$\int_{M^n}\varphi_{\alpha}u_{\beta}dv=0,\ \ \ \
1\leq\beta \leq  \alpha\leq N.
$$
Thus, 
from the Rayleigh-Ritz inequality, we have
\begin{equation}
\lambda_{\alpha+1}{\int_{M^n}\varphi_{\alpha}^2}dv\leq{\int_{M^n}\varphi_{\alpha}P_g\varphi_{\alpha}}dv,\
\ 1\leq\alpha\leq N.
\end{equation}
\begin{equation*}
\aligned
P_g\varphi_{\alpha}&= P_g(z_{\alpha}u_1) - a_{\alpha}P_gu_1=P_g(z_{\alpha}u_1) - \lambda_1a_{\alpha}u_1.
\endaligned
\end{equation*}

\begin{equation*}
\aligned
&P_g(z_{\alpha}u_1)\\
= &\Delta^2(z_{\alpha}u_1) -\text{div}\bigl [(a_nRg+b_n\text{Ric})\nabla (z_{\alpha}u_1)\bigl]+\dfrac{n-4}2Q(z_{\alpha}u_1)\\
=&\Delta^2z_{\alpha}\ u_1 +2\Delta z_{\alpha}\Delta u_1 +2\Delta g( \nabla z_{\alpha}, \nabla  u_1)\\
&+2g(\nabla z_{\alpha}, \nabla (\Delta u_1))+z_{\alpha}\Delta^2 u_1+2g(\nabla (\Delta z_{\alpha}), \nabla u_1)\\
&-\text{div}\bigl [u_1(a_nRg+b_n\text{Ric})\nabla z_{\alpha}\bigl]
-\text{div}\bigl [z_{\alpha}(a_nRg+b_n\text{Ric})\nabla u_1\bigl]+\dfrac{n-4}2Q(z_{\alpha}u_1)\\
=&\Delta^2z_{\alpha} \ u_1 +2\Delta z_{\alpha}\Delta u_1 +2\Delta g( \nabla z_{\alpha}, \nabla  u_1)
+2g(\nabla z_{\alpha}, \nabla (\Delta u_1))+2g(\nabla (\Delta z_{\alpha}), \nabla u_1)\\
&-\text{div}\bigl [u_1(a_nRg+b_n\text{Ric})\nabla z_{\alpha}\bigl]
-g(\nabla z_{\alpha}, (a_nRg+b_n\text{Ric})\nabla u_1)+z_{\alpha}Pu_1\\
=&r_{\alpha}+\lambda_1z_{\alpha}u_1
\endaligned
\end{equation*}
with
\begin{equation*}
\begin{aligned}
r_{\alpha}
&=\Delta^2z_{\alpha} \ u_1 +2\Delta z_{\alpha}\Delta u_1 +2\Delta g( \nabla z_{\alpha}, \nabla  u_1)\\
&+2g(\nabla z_{\alpha}, \nabla (\Delta u_1))+2g(\nabla (\Delta z_{\alpha}), \nabla u_1)\\
&-\text{div}\bigl [u_1(a_nRg+b_n\text{Ric})\nabla z_{\alpha}\bigl]
-g(\nabla z_{\alpha}, (a_nRg+b_n\text{Ric})\nabla u_1).\\
\end{aligned}
\end{equation*}
According to the Stokes formula, we derive
$$
\int_{M^n}r_{\alpha}u_1dv=0.
$$
Letting 
$$
w_{\alpha}=\int_{M^n}r_{\alpha}\varphi_{\alpha}dv
$$

\begin{equation*}
\aligned
\int_{M^n}\varphi_{\alpha}P_g\varphi_{\alpha}dv&
=\int_{M^n}\varphi_{\alpha}\bigl(P_g(z_{\alpha}u_1)-\lambda_1a_{\alpha}u_1\bigl)dv\\
&=\int_{M^n}\varphi_{\alpha}\bigl(r_{\alpha}+\lambda_1\varphi_{\alpha}\bigl)dv.\\
\endaligned
\end{equation*}
Hence,
\begin{equation}
(\lambda_{\alpha+1}-\lambda_1){\int_{M^n}\varphi_{\alpha}^2}dv\leq{\int_{M^n}\varphi_{\alpha}r_{\alpha}}dv
=w_{\alpha}=\int_{M^n}z_{\alpha}u_1r_{\alpha}dv,\
\ 1\leq\alpha\leq N.
\end{equation}
By a direct calculation, we obtain
\begin{equation*}
\begin{aligned}
&2\int_{M^n}z_{\alpha}u_1g(\nabla (\Delta z_{\alpha}), \nabla u_1)dv=\int_{M^n}(\Delta z_{\alpha})^2u_1^2dv\\
&+\int_{M^n}\Delta z_{\alpha}g(\nabla z_{\alpha},\nabla u_1^2)dv
-\int_{M^n}(z_{\alpha}\Delta^2 z_{\alpha})u_1^2dv,
\end{aligned}
\end{equation*}

\begin{equation*}
\begin{aligned}
&2\int_{M^n}z_{\alpha}u_1\Delta g( \nabla z_{\alpha}, \nabla  u_1)dv=2\int_{M^n}u_1\Delta z_{\alpha}g( \nabla z_{\alpha}, \nabla  u_1)dv\\
&+2\int_{M^n}z_{\alpha}\Delta u_1g( \nabla z_{\alpha}, \nabla  u_1)dv+4\int_{M^n}g(\nabla z_{\alpha},\nabla u_1)^2dv
\end{aligned}
\end{equation*}

\begin{equation*}
\begin{aligned}
&2\int_{M^n}z_{\alpha}u_1g(\nabla  z_{\alpha}, \nabla (\Delta u_1))dv=-2\int_{M^n}u_1z_{\alpha}\Delta z_{\alpha}\Delta u_1dv\\
&-2\int_{M^n}u_1\Delta u_1g(\nabla z_{\alpha},\nabla z_{\alpha})dv
-2\int_{M^n}z_{\alpha}g(\nabla z_{\alpha},\nabla u_1)\Delta u_1dv.
\end{aligned}
\end{equation*}
Thus, we derive
\begin{equation}
\begin{aligned}
&(\lambda_{\alpha+1}-\lambda_1){\int_{M^n}\varphi_{\alpha}^2}dv\leq 
w_{\alpha}=\int_{M^n}z_{\alpha}u_1r_{\alpha}dv\\
&=\int_{M^n}\bigl(u_1\Delta z_{\alpha}+2g( \nabla z_{\alpha}, \nabla  u_1)\bigl)^2dv+\int_{M^n}u_1^2g\bigl((a_nRg
+b_n\text{Ric})\nabla z_{\alpha},\nabla z_{\alpha}\bigl)dv\\
&-2\int_{M^n}g(\nabla z_{\alpha},\nabla z_{\alpha})u_1\Delta u_1dv,  \
\ 1\leq\alpha\leq N.
\end{aligned}
\end{equation}
From the lemma 2.1, we have
\begin{equation}
\aligned
&\sum_{\alpha=1}^N(\lambda_{\alpha+1}-\lambda_1){\int_{M^n}\varphi_{\alpha}^2}dv\\
&\leq \int_{M^n}\bigl(n^2|H|^2+(na_n+b_n)R\bigl)u_1^2dv+2(n+2)\int_{M^n}g( \nabla u_1, \nabla  u_1)dv. 
\endaligned
\end{equation}

On the other hand, 
\begin{equation}
\aligned &\int_{M^n}\varphi_{\alpha}\biggl(u_1\Delta z_{\alpha}+2g( \nabla z_{\alpha}, \nabla  u_1)\biggl)dv\\
&=\int_{M^4}( z_{\alpha}-a_{\alpha})u_1\biggl(u_1\Delta z_{\alpha}+2g( \nabla z_{\alpha}, \nabla  u_1)\biggl)dv\\
                     &=-\int_{M^4}|u_1\nabla z_{\alpha}|^2dv.
\endaligned
\end{equation}
Therefore, for any positive $\delta>0$, we obtain, from  (3.5), 
\begin{equation}
 \aligned &
 (\lambda_{\alpha+1}-\lambda_1)^{\frac12}\int_{M^n}|u_1\nabla z_{\alpha}|^2dv\\
&=-(\lambda_{\alpha+1}-\lambda_1)^{\frac12}\int_{M^n}
\varphi_{\alpha}\biggl (u_1\Delta z_{\alpha}+2g( \nabla z_{\alpha}, \nabla  u_1)\biggl)dv\\
 &\leq \frac12\biggl\{\delta(\lambda_{\alpha+1}-\lambda_1)\int_{M^n}\varphi_{\alpha}^2dv
 +\dfrac{1}{\delta} \int_{M^n}\biggl(u_1\Delta z_{\alpha}+2g( \nabla z_{\alpha}, \nabla  u_1)\biggl)^2dv\biggl\}.\\
 \endaligned
\end{equation}
According to (3.4) and (3.6), we infer
\begin{equation}
 \aligned &
 \sum_{\alpha=1}^N(\lambda_{\alpha+1}-\lambda_1)^{\frac12}\int_{M^n}|u_1\nabla z_{\alpha}|^2dv\\
 &\leq \frac12\biggl\{\delta\sum_{\alpha=1}^N(\lambda_{\alpha+1}-\lambda_1)\int_{M^n}\varphi_{\alpha}^2dv
 +\dfrac{1}{\delta} \sum_{\alpha=1}^N\int_{M^n}\bigl(u_1\Delta z_{\alpha}+2g( \nabla z_{\alpha}, \nabla  u_1)\bigl)^2dv\biggl\}\\
 &\leq \frac12\delta\biggl\{\int_{M^n}\bigl(n^2|H|^2+(na_n+b_n)R\bigl)u_1^2dv
 +2(n+2)\int_{M^n}g( \nabla u_1, \nabla  u_1)dv\biggl\}\\
 &+\dfrac{1}{2\delta}\biggl\{\int_{M^n}n^2|H|^2u_1^2dv+4\int_{M^n}g( \nabla u_1, \nabla  u_1)dv\biggl\}.\\
 \endaligned
\end{equation}
By the same proof as the formula (2.5) in the section 2, we have 
\begin{equation}
\aligned
&\sum_{\alpha=1}^N(\lambda_{\alpha+1}-\lambda_1)^{\frac{1}{2}}|\nabla z_{\alpha}|^2
\geq \sum\limits_{j=1}^n(\lambda_{j+1}-\lambda_1)^{\frac{1}{2}}.
\endaligned
\end{equation}
Hence, 
we obtain
\begin{equation}
 \aligned 
  &\sum\limits_{j=1}^n(\lambda_{j+1}-\lambda_1)^{\frac{1}{2}}\\
  &\leq \frac12\delta\biggl(\int_{M^n}\bigl(n^2|H|^2+(na_n+b_n)R\bigl)u_1^2dv+2(n+2)\int_{M^n}g( \nabla u_1, \nabla  u_1)dv\biggl)\\
 &+\dfrac{1}{2\delta}\biggl(\int_{M^n}n^2|H|^2u_1^2dv+4\int_{M^n}g( \nabla u_1, \nabla  u_1)dv\biggl).\\
 \endaligned
\end{equation}
Letting $S$ denote the squared norm of the second fundamental form of $M^n$, from the Gauss equation, we have 
$$
R=n(n-1)|H|^2-(S-n|H|^2)\leq n(n-1)|H|^2.
$$
Since 
$$
na_n+b_n=\dfrac{n^2-2n-4}{2(n-1)}>0,
$$
we have
$$
n^2|H|^2+(na_n+b_n)R\leq \dfrac{n(n^2-4)|H|^2}2.
$$
Taking 
$$
\dfrac1{\delta} =\sqrt{\dfrac{\int_{M^n}n^2|H|^2u_1^2dv+4\int_{M^n}g( \nabla u_1, \nabla  u_1)dv}
{ \int_{M^n}\dfrac{n(n^2-4)|H|^2}2u_1^2dv+2(n+2)\int_{M^n}g( \nabla u_1, \nabla  u_1)dv}}
 $$
 we have
 \begin{equation}
 \aligned 
  &\sum\limits_{j=1}^n(\lambda_{j+1}-\lambda_1)^{\frac{1}{2}}\\
  &\leq\sqrt{\int_{M^n}\dfrac{n(n^2-4)|H|^2}2u_1^2dv+2(n+2)\int_{M^n}g( \nabla u_1, \nabla  u_1)dv}\\
&\times\sqrt{\int_{M^n}n^2|H|^2u_1^2dv+4\int_{M^n}g( \nabla u_1, \nabla  u_1)dv}
\endaligned
\end{equation}
If the equality holds, we have 
$$
\lambda_2=\lambda_3=\cdots=\lambda_N,
$$
and $S\equiv n|H|^2$. Thus, $M^n$ is totally umbilical, that is, $M^n$ is isometric to a sphere.
It  completes the proof of the theorem 1.2.

\begin{corollary}
Let $(M^n, g)$ $(n>4)$ be an  $n$-dimensional compact  submanifold with the 
 metric $g$ in  the unit sphere $S^{N}(1)$. Then, eigenvalues 
of the Paneitz operator $P_g$ satisfy
 \begin{equation}
 \aligned 
  &\sum\limits_{j=1}^n(\lambda_{j+1}-\lambda_1)^{\frac{1}{2}}\\
  &\leq\sqrt{\int_{M^n}\dfrac{n(n^2-4)(|H|^2+1)}2u_1^2dv+2(n+2)\int_{M^n}g( \nabla u_1, \nabla  u_1)dv}\\
&\times\sqrt{\int_{M^n}n^2(|H|^2+1) u_1^2dv+4\int_{M^n}g( \nabla u_1, \nabla  u_1)dv}
\endaligned
\end{equation}
and  the equality holds if and only if  $M^n$ is isometric to a sphere $S^n(r)$,
where $u_1$ is the normalized first eigenfunction of $P_g$.
\end{corollary}

\vskip 3pt\noindent {\it Proof of  Theorem 1.3}. 
Since $n\neq 4$, we assume  that eigenvalues of the Paneitz operator $P_g$
satisfy
$$
0<\lambda_1< \lambda_2\leq \cdots, \lambda_k\leq \cdots\to +\infty.
$$
Let
$u_i$ be an eigenfunction corresponding to  eigenvalue
$\lambda_i$ such that $\{u_i\}_{i=1}^{\infty}$ becomes an
orthonormal basis of $L^2(M^n)$, that is,
\begin{equation*}
\begin{cases}
P_gu_i=\lambda_i u_i, \\
\int_{M^n}u_iu_jdv=\delta_{ij}, \ \ i, j= 1, 2,  \cdots.
\end{cases}
\end{equation*}
We shall use the similar  method to prove the theorem 1.3. 
We define an $(N+1)\times (N+1)$-matrix $A$ as follows:
$$
A:=(a_{\alpha\beta})
$$
where $a_{\alpha\beta}=\int_{M^n}y_{\alpha}u_{\beta}dv$, for $\alpha, \beta =1, 2, \cdots, N+1$,  and
${\bf y}=(y_{\alpha})$ is the position vector of the immersion in
$\mathbf{R}^{N+1}$ with $|{\bf y}|^2=\sum_{\alpha=1}^{N+1}y_{\alpha}^2=1$. 
Thus, there is  an orthogonal matrix  $U=(q_{\alpha\beta})$ such that 
$$
\int_{M^n}z_{\alpha}u_{\beta}dv=0,\ \ 1\leq\beta<\alpha\leq N+1,
$$
where 
$z_{\alpha}=\sum\limits_{\gamma=1}^{N+1}q_{\alpha\gamma}y_{\gamma}$.
Since $U$ is an orthogonal matrix, we have 
$$
\sum_{\alpha=1}^{N+1}z_{\alpha}^2=1.
$$
Putting 
$$
\psi_{\alpha}:=z_{\alpha},\ \ \  \ \ \ \ \ \ 1\leq
\alpha\leq N+1,
$$ 
we infer 
$$\int_{M^n}\psi_{\alpha}u_{\beta}dv=0,\ \ \ \
1\leq\beta<\alpha\leq N+1.
$$
Thus, 
from the Rayleigh-Ritz inequality, we have
$$
\lambda_{\alpha}{\int_{M^n}\psi_{\alpha}^2}dv\leq{\int_{M^n}\psi_{\alpha}P_g\psi_{\alpha}}dv,\
\ 1\leq\alpha\leq N+1.
$$
\begin{equation}
\aligned
P_g\psi_{\alpha}&= P_g(z_{\alpha}).
\endaligned
\end{equation}
According to the Stokes formula, we derive
\begin{equation*}
\aligned
\int_{M^n}\psi_{\alpha}P_g\psi_{\alpha}dv&=\int_{M^n}\biggl[\bigl(\Delta z_{\alpha}\bigl)^2
+g((a_nRg+b_n\text{Ric})\nabla z_{\alpha}, \nabla z_{\alpha})
 +\dfrac{n-4}2Q(z_{\alpha})^2
\biggl]dv\\
       \endaligned
\end{equation*}
From the lemma 2.1, we have
\begin{equation}
\aligned
&\sum_{\alpha=1}^{N+1}\int_{M^n}\psi_{\alpha}P_g\psi_{\alpha}dv\\
&
=\sum_{\alpha=1}^{N+1}\int_{M^n}\biggl[\bigl(\Delta z_{\alpha}\bigl)^2
+g((a_nRg+b_n\text{Ric})\nabla z_{\alpha}, \nabla z_{\alpha})
 +\dfrac{n-4}2Q(z_{\alpha})^2
\biggl]dv\\
&=\int_{M^n}\bigl((n^2|H|^2+n^2)+(na_n+b_n)R
 +\dfrac{n-4}2Q\bigl)dv.
       \endaligned
\end{equation}
Hence, 
\begin{equation}
\aligned
\sum_{\alpha=1}^{N+1}\lambda_{\alpha}{\int_{M^n}\psi_{\alpha}^2}dv
\leq \int_{M^n}\bigl((n^2|H|^2+n^2)+(na_n+b_n)R
 +\dfrac{n-4}2Q\bigl)dv. 
\endaligned
\end{equation}
On the other hand, 
\begin{equation}
\aligned &\int_{M^n}\psi_{\alpha}(\Delta z_{\alpha})dv=\int_{M^n}z_{\alpha}\Delta z_{\alpha}dv
=-\int_{M^n}|\nabla z_{\alpha}|^2dv.
\endaligned
\end{equation}
Therefore, for any positive $\delta>0$, we obtain
\begin{equation}
 \aligned &
 \lambda_{\alpha}^{\frac12}\int_{M^n}|\nabla z_{\alpha}|^2dv\\
&=-\lambda_{\alpha}^{\frac12}\int_{M^n}\psi_{\alpha}(\Delta z_{\alpha})dv\\
 &\leq \frac12\bigl(\delta\lambda_{\alpha}\int_{M^n}\psi_{\alpha}^2dv
 +\dfrac{1}{\delta} \int_{M^n}(\Delta z_{\alpha})^2dv\bigl)\\
 \endaligned
\end{equation}
and
\begin{equation}
 \aligned &
 \sum_{\alpha=1}^{N+1}\lambda_{\alpha}^{\frac12}\int_{M^n}|\nabla z_{\alpha}|^2dv\\
 &\leq \frac12\bigl(\delta\sum_{\alpha=1}^{N+1}\lambda_{\alpha}\int_{M^n}\psi_{\alpha}^2dv
 +\dfrac{1}{\delta}\sum_{\alpha=1}^{N+1} \int_{M^n}(\Delta z_{\alpha})^2dv\bigl)\\
 &\leq \frac12\biggl[\delta\int_{M^n}\bigl((n^2|H|^2+n^2)+(na_n+b_n)R
 +\dfrac{n-4}2Q\bigl)dv\\
 &+\dfrac{1}{\delta}\int_{M^n}(n^2|H|^2+n^2)dv\biggl].\\
 \endaligned
\end{equation}
By using the same proof as the formula (2.5) in the section 2, we have 
\begin{equation}
\aligned
&\sum_{\alpha=1}^{N+1}\lambda_{\alpha}^{\frac{1}{2}}|\nabla z_{\alpha}|^2
\geq\sum_{j=1}^n\lambda_{j}^{\frac{1}{2}}.
\endaligned
\end{equation}
Thus, we obtain
\begin{equation}
 \aligned 
  \sum_{j=1}^n\lambda_{j}^{\frac{1}{2}}\text{vol}(M^n)&\leq \frac12\biggl[\delta\int_{M^n}\bigl((n^2|H|^2+n^2)+(na_n+b_n)R
 +\dfrac{n-4}2Q\bigl)dv\\
 &+\dfrac{1}{\delta}\int_{M^n}(n^2|H|^2+n^2)dv\biggl].\\ 
 \endaligned
\end{equation}
Taking 
$$
\dfrac1{\delta} =\sqrt{\dfrac{\int_{M^n}\bigl((n^2|H|^2+n^2)+(na_n+b_n)R
 +\dfrac{n-4}2Q\bigl)dv}
{ \int_{M^n}(n^2|H|^2+n^2)dv}}
 $$
 we have
\begin{equation*}
 \aligned &
 \sum_{j=1}^n\lambda_{j}^{\frac{1}{2}}\leq n\dfrac{\sqrt{\int_{M^n}\bigl((n^2|H|^2+n^2)+(na_n+b_n)R
 +\dfrac{n-4}2Q\bigl)dv\int_{M^n}(|H|^2+1)dv}}{\text{vol}(M^n)}.\\
 \endaligned
\end{equation*}
If the equality holds, we have 
$$
\lambda_2=\lambda_3=\cdots=\lambda_{N+1},
$$
  $|\nabla z_1|\equiv1$ because of $\lambda_1<\lambda_2$ and 
  $$
\Delta z_{1}=-\sqrt{\lambda_{1}}\delta  z_{1}, \  \ 
\Delta z_{\alpha}=-\sqrt{\lambda_{n}}\delta  z_{\alpha} \ \  \text{for} \  \alpha>1.
$$
Since 
$$
\sum_{\alpha=1}^{N+1}z_{\alpha}^2=1, 
$$
we have 
$$
n-\sqrt{\lambda_{n}}\delta +(\sqrt{\lambda_{n}}\delta -\sqrt{\lambda_{1}}\delta) z_{1}^2=0.
$$
Thus, 
$$
\sqrt{\lambda_{n}}\delta=\sqrt{\lambda_{1}}\delta 
$$
or $z_1^2$ is constant.
It is impossible because  $|\nabla z_1|\equiv1$ and  $\lambda_1<\lambda_2$.
Therefore, the equality does not hold. It  completes the proof of the theorem 1.3.

\end{document}